\documentclass[a4paper, 11pt, oneside, onecolumn,sort&compress]{elsarticle}
\usepackage{amssymb}
\usepackage{mathrsfs}
\usepackage{amsfonts}
\usepackage{ntheorem}

\newtheorem*{thmA}{Theorem A}
\newtheorem*{thmB}{Theorem B}

\newtheorem{theorem}{Theorem}[section]

\newdefinition{example}{Example}[section]
\newdefinition{remark}{Remark}[section]
\newdefinition{definition}{Definition}[section]
\newproof{proof}{Proof}[section]
\newproof{pot}{Proof of Theorem \ref{the1.1}}
\usepackage{amssymb,amsmath}

\begin{document}
\begin{frontmatter}

\title{\textbf{Multiple solutions for Kirchhoff equations under the partially sublinear case}
}
\author  {Xiaojing Feng\corref{cor1}}
\ead{fengxj@sxu.edu.cn}

\cortext[cor1]{Corresponding author.}

\address{School of Mathematical Sciences, Shanxi University, Taiyuan 030006, People's Republic of China}

\begin{abstract}

In this paper, we prove the infinitely many solutions to a class of
sublinear Kirchhoff type equations by using an extension of Clark's
theorem established by Zhaoli Liu and Zhi-Qiang Wang.
\end{abstract}
\end{frontmatter}

\section{Introduction and main results}
In this paper we study the existence and multiplicity of solutions
for the following Kirchhoff type equations:
$$\left(a+\int_{\mathbb{R}^N}|\nabla u|^2+b\int_{\mathbb{R}^N}u^2\right)[-\Delta u+bu]
=K(x)f(x,u),\ {\rm in}\ \mathbb{R}^3,\eqno{(1.1)}$$
where $a,\ b$ are positive constants.

When $\Omega$ is a smooth bounded domain in $\mathbb{R}^3$, the problem
$$\left\{\begin{array}{ll}
-\left(a+b\int_{\Omega}|\nabla u|^2dx\right)\triangle u=f(x,u),\ &{\rm in}\ \Omega\\
u=0,&{\rm on}\ \partial\Omega,
\end{array}\right.\eqno{(1.2)}$$
has been many papers concerned.
Perera and Zhang \cite{pz1} considered the case
where $f(x,\cdot)$ is asymptotically linear at $0$ and asymptotically $4-$linear at infinity. They obtained a nontrivial solution of the
problems by using the Yang index and critical group.
Then, in \cite{pz1} they considered the cases where $f(x,\cdot)$ is $4-$sublinear,
$4-$superlinear and asymptotically $4-$linear at infinity. By various assumption on $f(x,\cdot)$ near $0$,
they obtained multiple and sign changing solutions.
Cheng and Wu \cite{cw}, Ma and Rivera \cite{mr} studied the existence of positive solutions of (1.2) and
He and Zou \cite{hz} obtained the existence of infinitely many positive solutions of (1.2), respectively;
Mao and Luan \cite{ml}
obtained the existence of signed and sign-changing solutions for the problem (1.2) with asymptotically $4-$linear bounded
nonlinearity via variational methods and invariant sets of descent flow; Sun and Tang \cite{st} studied the existence and multiplicity
results of nontrivial solutions for the problem (1.2) with the weaker monotony and $4-$superlinear nonlinearity.
For (1.2), Sun and Liu \cite{sl} considered the
cases where the nonlinearity is superlinear near zero but asymptotically $4-$linear at infinity,
and the nonlinearity is asymptotically linear near zero but $4-$superlinear at infinity. By
computing the relevant critical groups, they obtained nontrivial solutions via Morse theory.

Comparing with (1.1) and (1.2), $\mathbb{R}^3$ in place of the bounded domain $\Omega\subset \mathbb{R}^3$.
This make that the study of the problem (1.1) is more difficult and interesting.
Wu \cite{xw} considered a class of Schr\"{o}dinger Kirchhoff-type problem in $\mathbb{R}^N$ and a sequence of high
energy solutions are obtained by using a symmetric Mountain Pass
Theorem. In \cite{af}, Alves and Figueiredo study a periodic Kirchhoff equation in $\mathbb{R}^N$, they get the nontrivial solution when
the nonlinearity is in subcritical case and critical case. Liu and He \cite{lh} get multiplicity of high energy solutions
for superlinear Kirchhoff equations in $\mathbb{R}^3$. Li, Li and Shi in \cite{ll} proved the existence of a positive solution
to a Kirchhoff type problem on $\mathbb{R}^N$ by using variational methods and cut-off functional technique.

In \cite{jw}, Jin and Wu in consider the following problem:
$$\left\{\begin{array}{ll}
-\left(a+b\int_{\mathbb{R}^N}|\nabla u|^2dx\right)\triangle u+u=f(x,u),\ &{\rm in}\ \mathbb{R}^N,\\
u\in H^1(\mathbb{R}^N),
\end{array}\right.\eqno{(1.3)}$$
where constants $a>0, b>0$, $N=2$ or $3$ and $f\in C(R^N\times R, R)$.

By using the
Fountain Theorem, they obtained the following theorem.

\begin{thmA}\cite{jw} Assume that the following conditions hold:

If the following assumptions are satisfied,

(H$_1$) $f(x,u)=o(|u|)$ as $|u|\to 0$ uniformly for any $x\in \mathbb{R}^N$.

(H$_2$) There are constants $1<p<2^*-1$ and $c>0$ such that
$$|f(x,u)|\leq c(1+|u|^p),\ \forall (x,u)\in \mathbb{R}^N\times \mathbb{R},$$
where
$$2^*-1=\left\{\begin{array}{ll}
\frac{N+2}{N-2},&N\geq 3;\\+\infty,&N=1,2.
\end{array}\right.$$

(H$_3$) There exists $\mu>4$ such that
$$\mu F(x,u)=\mu\int_0^uf(x,s)ds\leq uf(x,u),\ \forall (x,u)\in \mathbb{R}^N\times \mathbb{R}.$$

(H$_4$) $$\inf_{x\in \mathbb{R}^N,|u|=1}F(x,u)>0$$

(H$_5$) $f(gx,u)=f(x,u)$ for each $g\in O(N)$ and for each $(x,u)\in \mathbb{R}^N\times \mathbb{R}$, where
$O(N)$ is the group of orthogonal transformations on $\mathbb{R}^N$.

(H$_6$) $f(x,-u)=-f(x,u)$ for any $(x,u)\in \mathbb{R}^N\times \mathbb{R}$.

Then problem (1.3) has a sequence $\{u_k\}$ of radial solutions.

\end{thmA}

Recently, the authors obtained an extension of Clark's theorem as follows.

\begin{thmB}\cite{lw}
Let $X$ be a Banach space, $\Phi\in C^1(X,\mathbb{R})$. Assume $\Phi$ is even and satisfies the (PS) condition,
bounded from below, and $\Phi(0)=0$. If for any $k\in \mathbb{N}$, there exists a $k-$dimensional subspace $X^k$
of $X$ and $\rho_k>0$ such that $\sup_{X^k\cap S_{\rho_k}}\Phi<0$, where $S_\rho=\{u\in X|\|u\|=\rho\}$, then
at least one of the following conclusions holds.

(i) There exists a sequence of critical points $\{u_k\}$ satisfying $\Phi(u_k)<0$ for all $k$ and
$\|u_k\|\to 0$ as $k\to \infty$.

(ii) There exists $r>0$ such that for any $0<a<r$ there exists a critical point $u$ such that $\|u\|=a$
and $\Phi(u)=0$.
\end{thmB}

In this paper, we consider the multiple solutions for Kirchhoff equations under the partially sublinear case
by using the Theorem C. Our main result is as follows.
\begin{theorem}\label{the1.1}
Assume that $f$ satisfies ($B_3$) and the following conditions:

($f_1$) There exist $\delta>0,\ 1\leq \gamma<2,\ C>0$ such that $f\in
C(\mathbb{R}^3\times[-\delta,\delta],\mathbb{R})$ and
$|f(x,z)|\leq C|z|^{\gamma-1}$;

($f_2$) $\lim_{z\to 0}F(x,z)/|z|^2=+\infty$
uniformly in some ball $B_r(x_0)\subset\mathbb{R}^3$, where $F(x,z)=\int_0^zf(x,s)ds$.

($f_3$) $K:\mathbb{R}^3\to \mathbb{R}^+$
is a positive continuous function such that $K\in L^{2/(1-\gamma)}(\mathbb{R}^3)\cap L^{\infty}(\mathbb{R}^3)$.

Then (1.1) possesses infinitely many solutions $\{u_k\}$ such that
$\|u_k\|_{L^\infty}\to 0$ as $k\to \infty$.
\end{theorem}

\begin{remark} Throughout the paper we denote by $C>0$ various positive constants which may
vary from line to line and are not essential to the problem.
\end{remark}
The paper is organized as follows: in Section 2, some preliminary results are presented.
Section 3 is devoted to the proof of Theorem 1.1.

\section{Preliminary\label{Section 2}}

In this Section, we will give some notations and Lemma that will be used throughout this paper.

Let $H^1=H^1(\mathbb{R}^3)$ be the completion of $C_0^\infty(\mathbb{R}^3)$ with respect to the inner product and norm
$$(u,v)=\int_{\mathbb{R}^3}[\nabla u\nabla v+buvdx,\,\ \|u\|=(u,u)^{1/2}.$$
Moreover, we denote the completion of $C_0^\infty(\mathbb{R}^3)$ with respect to the norm
$$\|u\|_{D^1}=\int_{\mathbb{R}^3}|\nabla u|^2dx$$
by $D^1=D^1(\mathbb{R}^3)$.
To avoid lack of compactness, we need consider the set of radial functions as follows:
$$H=H_r^1(\mathbb{R}^3)=\{u\in H^1(\mathbb{R}^3)| u(x)=u(|x|)\}.$$
Here we note that the continuous embedding $H\hookrightarrow L^q(\mathbb{R}^3)$
is compact for any $q\in (2,6)$.

Define a functional $$J_1(u)=\frac{1}{2}\|u\|^2+\frac{1}{4}\|u\|^4-\int_{\mathbb{R}^3}K(x)F(x,u),\ u\in H.$$
Then we have from (f$_1$) that $J_1$ is well defined on $H$ and is of $C^1$, and
$$(J_1(u),v)=a(u,v)+\|u\|^2(u,v)-\int_{\mathbb{R}^3}K(x)f(x,u)v,\ u,\ v\in H.$$
It is standard to verify that the weak solutions of (1.1) correspond to the critical points of the functional $J_1$.

\section{Proofs of the main result \label{Section 3}}

{\bf Proof of Theorem 1.1.} Choose $\hat{f}\in C(\mathbb{R}^N\times\mathbb{R}, \mathbb{R})$ such that $\hat{f}$ is
odd in $u\in \mathbb{R}$, $\hat{f}(x,u)=f(x,u)$ for $x\in \mathbb{R}^N$ and $|u|<\delta/2$, and
$\hat{f}(x,u)=0$ for $x\in \mathbb{R}^N$ and $|u|>\delta$. In order to obtain solutions of (1.1) we consider

$$\left(a+\int_{\mathbb{R}^N}|\nabla u|^2+b\int_{\mathbb{R}^N}u^2\right)[-\Delta u+bu]
=K(x)\hat{f}(x,u),\ {\rm in}\ \mathbb{R}^N,\eqno{(1.1)}$$

Moreover, (3.1) is variational and its solutions are the critical
points of the functional defined in $H$ by
$$J(u)=\frac{1}{2}a\|u\|^2+\frac{1}{4}\|u\|^4-\int_{\mathbb{R}^3}K(x)\hat{F}(x,u)dx.$$
From ($f_1$), it is easy to check that $J$ is well defined on $H$ and
$J\in C^1(H^1(\mathbb{R}^3),\mathbb{R})$ (see \cite{ct} for more detail), and
$$J'(u)v=a(u,v)+\|u\|^2(u,v)
-\int_{\mathbb{R}^3}K(x)\hat{f}(x,u)vdx,\ v\in H.$$

Note that $J$ is even, and $J(0)=0$.
For $u\in H^1(\mathbb{R}^3)$,
$$\int_{\mathbb{R}^3}K(x)|\hat{F}(x,u)|dx\leq C\int_{\mathbb{R}^3}K(x)|u|^\gamma dx
\leq
C\|K\|_{L^{\frac{2}{2-\gamma}}(\mathbb{R}^3)}\|u\|^\gamma_{L^{2}(\mathbb{R}^3)}\leq
C\|u\|^\gamma.$$ Hence, it follows from Lemma 2.1 that
$$J(u)\geq\frac{1}{2}\|u\|^2-C\|u\|^\gamma,\ u\in H.\eqno(3.2)$$
We now use the same ideas to prove the (PS) condition. Let $\{u_n\}$ be a sequence in $H$ so that
$J(u_n)$ is bounded and $J'(u_n)\to 0$. We shall prove that $\{u_n\}$ converges. By (3.2), we claim that $\{u_n\}$
is bounded.
Assume without loss of generality that $\{u_n\}$ converges to $u$ weakly in $H$. Observe that
\begin{equation*}
\begin{split}
\langle J'(u_n)-J'(u),u_n-u\rangle&=a\|u_n-u\|^2+\|u_n\|^2\|u_n-u\|^2\\
&+(\|u_n\|^2-\|u\|^2)(u,u_n-u)\\
&-\int_{\mathbb{R}^3}K(x)(\hat{f}(x,u_n)-\hat{f}(x,u))(u_n-u)dx.\\
\end{split}
\end{equation*}
Hence, we have
\begin{equation*}
\begin{split}
a\|u_n-u\|^2&\leq\langle J'(u_n)-J'(u),u_n-u\rangle-(\|u_n\|^2-\|u\|^2)(u,u_n-u)\\
&+\int_{\mathbb{R}^3}K(x)(\hat{f}(x,u_n)-\hat{f}(x,u))(u_n-u)dx.\\
&\equiv I_1+I_2+I_3,
\end{split}
\end{equation*}

It is clear that $I_1\to 0$ and $I_2\to 0$ as $n\to \infty.$ In the following, we will estimate $I_3$, by using
(f$_3$),for any $R>0$,
\begin{equation*}
\begin{split}
&\int_{\mathbb{R}^3}K(x)|\hat{f}(x,u_n)-\hat{f}(x,u)||u_n-u|dx\\
&\leq C\int_{\mathbb{R}^3\setminus B_R(0)}K(x)(|u_n|^\gamma+|u|^\gamma)dx+C\int_{B_R(0)}(|u_n|^{\gamma-1}+|u|^{\gamma-1})|u_n-u|dx\\
&\leq C\left(\|u_n\|^\gamma_{L^2(\mathbb{R}^3\setminus B_R(0))}+\|u\|^\gamma_{L^2(\mathbb{R}^3\setminus B_R(0))}\right)
\|K\|_{L^{\frac{2}{2-\gamma}}(\mathbb{R}^3\setminus B_R(0))}\\
&+C\left(\|u_n\|^{\gamma-1}_{L^{\gamma}(B_R(0))}+\|u\|^{\gamma-1}_{L^{\gamma}(B_R(0))}\right)\|u_n-u\|_{L^{\gamma}(B_R(0))}\\
&\leq C\|K\|_{L^{\frac{2}{2-\gamma}}(\mathbb{R}^3\setminus B_R(0))}+C\|u_n-u\|_{L^{\gamma}(B_R(0))},
\end{split}
\end{equation*}
which implies
$$\lim_{n\to+\infty}\int_{\mathbb{R}^3}K(x)|\hat{f}(x,u_n)-\hat{f}(x,u))||u_n-u|dx=0.$$
Therefore, $\{u_n\}$ converges strongly in $H$ and the (PS) condition holds for $J$. By (f$_2$) and (f$_3$),
for any $L>0$, there exists
$\delta=\delta(L)>0$ such that if $u\in C_0^\infty(B_r(x_0))$ and $|u|_\infty<\delta$
then $K(x)\hat{F}(x,u(x))\geq L|u(x)|^2$,
and it follows from Lemma 2.1 that
$$J(u)\leq\frac{a}{2}\|u\|^2+\frac{1}{4}\|u\|^4-L\|u\|^2_{L^2(\mathbb{R}^3)}.$$
This implies, for any $k\in \mathbb{N}$, if $X^k$ is a $k-$dimensional subspace of $C_0^\infty(B_r(x_0))$
and $\rho_k$ is sufficiently small then $\sup_{X^k\cap S_{\rho_k}}J(u)<0$, where $S_\rho=\{u\in \mathbb{R}^3| \|u\|=\rho\}$.
Now we apply Theorem C to obtain infinitely many solutions $\{u_k\}$ for (3.1) such that

$$\|u_k\|\to 0,\ k\to \infty.\eqno(3.3)$$

Finally we show that $\|u_k\|_{L^\infty}\to 0$ as $k\to \infty$. Let $u$ be a solution of (3.1) and
$\alpha>0$. Let $M>0$ and set $u^M(x)=max\{-M,min\{u(x),M\}\}$. Multiplying both sides of (3.1) with
$|u^M|^\alpha u^M$ implies
$$\frac{4a}{(\alpha+2)^2}\int_{\mathbb{R}^3}|\nabla|u^M|^{\frac{\alpha}{2}+1}|^2dx
\leq C\int_{\mathbb{R}^3}|u^M|^{\alpha+1}dx.$$
By using the iterating method in \cite{lw}, we can get the following estimate
$$\|u\|_{L^\infty(\mathbb{R}^3)}\leq C_1\|u\|^{\nu}_{L^6(\mathbb{R}^3)},$$
where $\nu$ is a number in $(0,1)$ and $C_1>0$ is independent of $u$ and $\alpha$. By (3.3) and
Sobolev imbedding Theorem\cite{w},
we derive that $\|u_k\|_{L^\infty(\mathbb{R}^3)}\to 0$ as $k\to \infty$. Therefore, $u_k$ are the solutions of (1.1) as
$k$ sufficiently large. The proof is completed.    \hfill$\Box$

\section{Acknowledgement \label{Section 4}}

 Authors would like to express their sincere gratitude to one anonymous referee for his/her
constructive comments for improving the quality of this paper.

\end{document}